\documentclass[12pt, a4paper,twoside]{article}
\usepackage{a4,amssymb,amsfonts,verbatim,pstricks,pst-plot,xypic,times}
\usepackage{yfonts}

\usepackage[english]{babel}
\usepackage{amssymb,latexsym}
\usepackage{amssymb}
\usepackage{latexsym}

\newcommand{\CC}{{\mathbb C}}
\newcommand{\RR}{{\mathbb R}}
\newcommand{\ZZ}{{\mathbb Z}}
\newcommand{\SSS}{{\mathbb S}}

\newcommand{\cZ}{{\mathcal Z}}

\newcommand{\n}{\nabla}

\renewcommand{\phi}{\varphi}

\newcommand{\eps}{\varepsilon}
\newcommand{\half}{\frac{1}{2}}
\renewcommand{\div}{\mathrm{div}}
\newcommand{\divM}{\mathrm{div}^M}
\newcommand{\divMt}{\mathrm{div}^{M_t}}
\newcommand{\grad}{\mathrm{grad}}
\newcommand{\gradM}{\mathrm{grad}^M}

\renewcommand{\Re}{\mathrm{Re}}
\newcommand{\End}{\mathrm{End}}
\newcommand{\scal}{\mathrm{Scal}}
\newcommand{\scalZ}{\mathrm{Scal}^\cZ}
\newcommand{\scalM}{\mathrm{Scal}^M}

\newcommand{\tr}{\mathrm{tr}}
\newcommand{\<}{\left\langle}       
\renewcommand{\>}{\right\rangle}       
\newcommand{\pa}[2]{\frac{\partial #1}{\partial #2}}
\newcommand{\SZ}{{\Sigma \cZ}}
\newcommand{\SZp}{{\Sigma^+ \cZ}}
\newcommand{\SZm}{{\Sigma^- \cZ}}
\newcommand{\SZpm}{{\Sigma^\pm \cZ}}

\newcommand{\ScN}{{\Sigma N}}

\newcommand{\SpcN}{{\Sigma_p N}}
\newcommand{\SM}{{\Sigma M}}

\newcommand{\nZ}{\nabla^\cZ}

\newcommand{\nM}{\nabla^M}
\newcommand{\nMt}{\nabla^{M_t}}

\newcommand{\nN}{\nabla^N}

\newcommand{\oN}{\omega^{N}}

\newcommand{\ncSN}{\nabla^{\Sigma N}}

\newcommand{\RN}{R^{N}}

\newcommand{\RSNc}{\mathcal{R}^{\Sigma N}}
\newcommand{\Spinc}{\mathrm{Spin^c}}

\newcommand{\RZ}{R^{\cZ}}

\newcommand{\RSZc}{\mathcal{R}^{\Sigma \cZ}}

\newcommand{\RicZ}{\mathrm{Ric}^\cZ}
\newcommand{\ricZ}{{\mathrm{ric}}^\cZ}

\newcommand{\RicN}{\mathrm{Ric}^N}

\newcommand{\ricN}{\mathrm{ric}^N}

\newcommand{\DDM}{\widetilde D}
\newcommand{\PSOM}{P_{\mathrm{SO}}M}

\newcommand{\PSON}{P_{\mathrm{SO}}N}
\newcommand{\PSOZ}{P_{\mathrm{SO}}\cZ}

\newcommand{\PSpincM}{P_{\mathrm{Spin^{c}}}M}

\newcommand{\PSpincN}{P_{\mathrm{Spin^{c}}}N}
\newcommand{\PSpincZ}{P_{\mathrm{Spin^{c}}}\cZ}
\newcommand{\PSZ}{P_{\SSS^1}\cZ}
\newcommand{\PSM}{P_{\SSS^1}M}

\newcommand{\SOrs}{\mathrm{SO}(r,s)} 

\newcommand{\SOrss}{\mathrm{SO}(r+1,s)}

\newcommand{\Clrs}{\mathrm{Cl}_{r,s}}
\newcommand{\CClrs}{\mathrm{\CC l}_{r,s}}

\newcommand{\CCllrs}{\mathrm{\CC l}_{r,s}^1}
 
\newcommand{\CClors}{\mathrm{\CC l}^0_{r,s}}

\newcommand{\Clrsp}{\mathrm{Cl}^*_{r,s}}

\newcommand{\CClorls}{\mathrm{\CC l}^0_{r+1,s}}
\newcommand{\Spinrs}{\mathrm{Spin}(r,s)} 
\newcommand{\Spinrsc}{\mathrm{Spin}^{c}(r,s)} 
\newcommand{\Spinrscc}{\mathrm{Spin}^{c}(r+1,s)} 

\newcommand{\Ad}{\mathrm{Ad}}
\newcommand{\Srs}{\Sigma_{r,s}}

\newcommand{\Sprs}{\Sigma_{r,s}^+}
\newcommand{\Smrs}{\Sigma_{r,s}^-}
\newcommand{\Sors}{\Sigma_{r,s}^0}

\newcommand{\Aut}{\mathrm{Aut}} 
\newcommand{\id}{\mathrm{Id}} 
 
\newcommand{\DW}{\mathfrak{D}^W}
\newcommand{\summe}{\sum_{j=1}^n}
\newcommand{\gdt}{\dot{g}_t}
\newcommand{\gddt}{\ddot{g}_t}

\newcommand{\cdott}{\bullet_t}
\newcommand{\Mt}{{M_t}}
\newcommand{\Mto}{{M_{t_0}}}

\newcommand{\DDMt}{{D}^{\Mt}}
\newcommand{\gradMt}{\mathrm{grad}^{\Mt}}
\newcommand{\DWt}{\mathfrak{D}^{W_t}}
\newcommand{\Dgdt}{\mathfrak{D}^{\gdt}}

\newtheorem{thm}{Theorem}[section]
\newtheorem{lemma}[thm]{Lemma}
\newtheorem{prop}[thm]{Proposition}
\newtheorem{cor}[thm]{Corollary}
\newtheorem{remark}[thm]{Remark}
\newtheorem{remarks}[thm]{Remarks}
\newtheorem{definition}[thm]{Definition}
\newtheorem{notation}[thm]{Notation}
\newtheorem{example}[thm]{Example}

\begin{document}
\title{The Energy-Momentum Tensor on $\Spinc$ Manifolds}
\author{Roger Nakad}
\maketitle
\begin{center}
Institut \'Elie Cartan, Universit\'e Henri Poincar\'e, Nancy I, B.P 239\\
54506 Vandoeuvre-L\`es-Nancy Cedex, France.
\end{center}
\begin{center}
 {\bf nakad@iecn.u-nancy.fr}
\end{center}

\begin{center}
 {\bf Abstract}
\end{center}
On $\Spinc$ manifolds, we study the Energy-Momentum tensor associated with a spinor field. First, we give a spinorial Gauss type formula for oriented hypersurfaces of a $\Spinc$ manifold. Using the notion of generalized cylinders, we derive the variationnal formula for the Dirac operator under metric deformation and point out that the Energy-Momentum tensor appears naturally as the second fundamental form of an isometric immersion. Finally, we show that generalized $\Spinc$ Killing spinors for Codazzi Energy-Momentum tensor are restrictions of parallel spinors.
\section{Introduction}
In \cite{hijazi}, O. Hijazi proved that on a compact Riemannian spin manifold $(M^n, g)$ any eigenvalue $\lambda$ of the Dirac operator to which is attached an eigenspinor $\psi$ satisfies
\begin{eqnarray}
 \lambda^2 \geqslant \inf_M (\frac 14 \scal^M + \vert \ell^\psi\vert^2),
\label{oussamascal}
\end{eqnarray}
where $\scal^M$ is the scalar curvature of the manifold $M$ and $\ell^\psi$ is the field of symmetric endomorphisms associated with the field of  quadratic forms $T^\psi$ called the Energy-Momentum tensor. It is defined on the complement set of zeroes of the eigenspinor $\psi$, for any vector $X \in \Gamma(TM)$ by 
$$T^\psi (X)= \Re < X\cdot\nabla_X \psi,\frac{\psi}{\vert\psi\vert^2}>.$$
Here $\nabla$ denotes the Levi-Civita connection on the spinor bundle of $M$ and ``$\cdot$'' the Clifford multiplication. The limiting case of (\ref{oussamascal}) is characterized by the existence of a spinor field $\psi$ satisfying for all $X \in \Gamma(TM)$,
\begin{eqnarray}
 \nabla_X\psi = -\ell^\psi(X)\cdot\psi.
\label{huhuhu}
\end{eqnarray}
For $\Spinc$ structures, the complex line bundle $L^M$ is endowed with an arbitrary connection and hence an arbitrary curvature $i\Omega^M$ which is an imaginairy $2$-form on the manifold. 
In terms of the Energy-Momentum tensor the author proved in \cite{article1} that on a compact Riemannian $\Spinc$ manifold any eigenvalue $\lambda$ of the  Dirac operator to which is attached an eigenspinor $\psi$ satisfies
\begin{eqnarray}
 \lambda ^2 \geqslant \inf_M\ \Big(\frac 14 \scal^M -\frac{c_n}{4}\vert \Omega^M\vert + \vert \ell^{\psi}\vert^2\Big),
\label{rodgy}
\end{eqnarray}
where $c_n =2[\frac n2]^\frac 12$. The limiting case of (\ref{rodgy}) is characterized by the existence of a spinor field $\psi$ satisfying for every $X \in \Gamma(TM)$,
\begin{eqnarray}
 \left\{
\begin{array}{l}
\n^\SM_X\psi=-\ell^\psi(X)\cdot\psi, \\ \\
\Omega^M \cdot\psi = i\frac{c_n}{2}\vert\Omega^M\vert\psi.
\end{array}
\right.
\end{eqnarray}
Here $\n^\SM$ denotes the Levi-Civita connection on the $\Spinc$ spinor bundle and ``$\cdot$'' the $\Spinc$ Clifford multiplication. In \cite{article1}, the author showed also that the sphere with a special $\Spinc$ structure is a limiting manifold for (\ref{rodgy}).\\ \\
Studying the Energy-Momentum tensor on a Riemannian or semi-Riemannian spin manifolds has been done by many authors, since it is related to several geometric constructions (see \cite{habib1}, \cite{3}, \cite{24} and \cite{friedrich} for results in this topic). In this paper we study the Energy-Momentum tensor on Riemannian and semi-Riemannian $\Spinc$ manifolds. First, we prove that the Energy-Momentum tensor appears in the study of the variations of the spectrum of the Dirac operator:
\begin{prop}
Let $(M^n,g)$ be a $\Spinc$ Riemannian manifold and $g_t = g+ tk$ a smooth 1-parameter family of metrics. For any spinor field $\psi\in \Gamma(\Sigma M)$, we have 
\begin{equation}
 \left.\frac{d}{dt}\right|_{t=0} (\DDMt \tau_{0}^t \psi, \tau_{0}^t \psi)_{g_t}
=-\frac 12\int_M < k, T_\psi > dv_g,
\end{equation}
where $( ., .) =\int_M \Re \< ., .\> dv_g$, the Dirac operator $\DDMt$ is the Dirac operator associated with $M_t = (M, g_t)$, $T_{\psi} =\vert\psi\vert^2\ T^\psi = \Re < X\cdot\nabla_X \psi,\psi>$ and $\tau_{0}^t \psi$ is the image of $\psi$ under the isometry  $\tau_{0}^t$ between the spinor bundles of $(M, g)$ and $(M, g_t)$.
\label{endomo}
\end{prop}
This was proven in \cite{1} by J. P. Bourguignon and P. Gauduchon for spin manifolds. Using this, we extend to $\Spinc$ manifolds a result by Th. Friedrich and E. C. Kim in \cite{4} on spin manifolds:
\begin{thm}
Let $M$ be a $\Spinc$ Riemannian manifold. A pair $(g_0,\psi_0)$ is a critical point of the Lagrange functional
$$\mathcal W(g,\psi)  = \int_U \big( \scal_g^M +\eps {\lambda <\psi,\psi>_g - < D_g\psi,\psi>_g}\Big)dv_g,$$
$(\lambda, \eps \in \RR)$ for all open subsets $U$ of $M$ if and only if $(g_0,\psi_0)$ is a solution of the following system
$$\left\{
\begin{array}{l}
\ \ \ \ \ \ \ \ \ \ \ \ \ \ \ \ \ \ D_g \psi = \lambda \psi,\\
\mathrm{ric}^M_g -\frac {\scal_g^M}{2}   \ g = \frac{\eps}{4}T_\psi,
\end{array}
\right.
$$
where $\mathrm{ric}^M_g$ denotes the Ricci curvature of $M$ considered as a symmetric bilinear form.
\label{frkim}
\end{thm}
Now, we interprete the Energy-Momentum tensor as the second fundamental form of a hypersurface. In fact, we prove the following:
\begin{prop}
Let $M^n \hookrightarrow (\cZ, g)$ be any compact oriented hypersurface isometrically immersed in an oriented Riemannian $\Spinc$ manifold $(\cZ, g)$, of constant mean curvature $H$ and Weingarten map $W$. Assume that $\cZ$ admits a parallel spinor field $\psi$, then the Energy-Momentum tensor associated with  $\phi =: \psi_{|_M}$ satisfies
$$2\ell^\phi =- W.$$
Moreover the hypersurface $M$ satisfies the equality case in (\ref{rodgy}) if and only if
\begin{eqnarray}
\scalZ - 2\ \ricZ(\nu, \nu)- c_n \vert\Omega^M\vert =0.
\label{maga}
\end{eqnarray}
\label{morel}
\end{prop}
\vskip-1cm
This was proven by Morel in \cite{24} for a compact oriented hypersurface of a spin manifold carrying parallel spinor but in this case the hypersurface $M$ is directly a limiting manifold for (\ref{oussamascal}) without the condition (\ref{maga}).
Finally, we study generalized Killing spinors on $\Spinc$ manifolds. They are characterized by the identity, for any tangent vector field $X$ on $M$,
\begin{eqnarray}
\n^\SM_X\psi=\frac 12 F(X)\cdot\psi,
\label{gks}
\end{eqnarray}
where $F$ is a given symmetric endomorphism on the tangent bundle. It is straightforward to see that 
$$ 2 T^\psi(X,Y) = -\< F(X), Y\>.$$
These spinors are closely related to the so-called $T$--Killing spinors studied 
by Friedrich and Kim in \cite{friedrich-kim99a} on spin manifolds. It is natural to ask whether the tensor $F$ can 
be realized as the Weingarten tensor of some isometric embedding of 
$M$ in a manifold $\cZ^{n+1}$ carrying parallel spinors. Morel studied this problem in the case of spin manifolds where the tensor $F$ is parallel and in \cite{3}, the authors studied the problem in the case of semi-Riemannian spin manifolds where the tensor $F$ is a Codazzi-Mainardi tensor. We establish the corresponding result for semi-Riemannian $\Spinc$ manifolds:
\begin{thm}\label{fini}
Let $(M^n,g)$ be a semi-Riemannian $\Spinc$ manifold carrying a generalized $\Spinc$ Killing spinor $\phi$ with a Codazzi-Mainardi tensor $F$. Then the generalized cylinder $\cZ :=I\times M$
with the metric $dt^2+g_t$, where $g_t(X,Y)=g((\id-tF)^2X,Y)$, equipped with the $\Spinc$ structure arising from the given one on $M$ has a parallel spinor whose restriction to $M$ is just $\phi$.
\end{thm}
A characterisation of limiting $3$-dimensional manifolds for (\ref{rodgy}), having generalized $\Spinc$ Killing spinors with Codazzi tensor is then given.\\ \\
The paper is organised as follows: In Section \ref{pre}, we collect basic material on spinors and the Dirac operator on semi-Riemannian $\Spinc$ manifolds. In Section \ref{hyper}, we study hypersurfaces of $\Spinc$ manifolds. We derive a spinorial Gauss formula after identifying the restriction of the $\Spinc$ spinor bundle of the ambient manifold with the $\Spinc$ spinor bundle of the hypersurface.  In Section \ref{gc}, we define the generalized cylinder of a $\Spinc$ manifold $M$  and  we collect formulas relating the curvature of a generalized cylinder to geometric data on $M$. In section \ref{vf}, we compare the Dirac operators for two differents semi-Riemannian metrics, then one first has to identify the spinor bundles using parallel transport. In the last section, we interprete the Energy-Momentum tensor as the second fundamental form of a hypersurface and we study generalized $\Spinc$ Killing spinors.
The author would like to thank Oussama Hijazi for his support and encouragements.
\section{The Dirac operator on semi-Riemannian Spin$^c$ manifolds}
\label{pre}
In this section, we collect some algebraic and geometric preliminaries concerning the Dirac operator on semi-Riemannian $\Spinc$ manifolds. Details can be found in \cite{baum} and \cite{3}. Let $r+s=n$ and consider on $\RR^n$ the nondegenerate symmetric bilinear form of signature $(r,s)$ given by
$$\< v,w \> := \sum_{j=1}^r v_j w_j - \sum_{j=r+1}^n v_j w_j,$$
for any $v$, $w \in \RR^n$. We denote by $\Clrs$ the real Clifford algebra corresponding to $(\RR^n, \<\cdot,\cdot\>)$, this is the unitary algebra generated by $\RR^n$ subject to the relations
$$e_j\cdot e_k + e_k\cdot e_j =
 \left\{
\begin{array}{l}
-2 \delta_{jk}\text{\ \ if\ \ }\ j\leqslant r, \\
 \ \ 2 \delta_{jk}\text{\ \ \ if\ \ }\ j > r,
\end{array}
\right.
$$
where $(e_j)_{1\leqslant j\leqslant n}$ is an orthonormal basis of $\RR^n$ of signature $(r,s)$, i.e., $\<e_j, e_k\> =\eps_j\delta_{jk}$ and $\eps_j =\pm 1$. The complex Clifford algebra  $\CClrs$ is the complexification of $\Clrs$ and it decomposes into even and odd elements $\CClrs = \CClors \oplus \CCllrs.$
The real spin group is defined by
$$
\Spinrs := \{ v_1\cdot...\cdot v_{2k} \in \Clrs\ |\ v_j \in \RR^n 
\mbox{ such that }\<v_j,v_j\> = \pm 1\}.
$$
The spin group $\Spinrs$ is the double cover of $\SOrs$, in fact the following  sequence is exact
$$1 \longrightarrow \ZZ/2\ZZ \longrightarrow \Spinrs
\stackrel{\xi}{\longrightarrow}  \SOrs  \longrightarrow 1,
$$
where $\xi= \Ad_{\arrowvert_{\Spinrs}}$ and $\Ad$ is defined by 
\begin{eqnarray*}
\Ad :\Clrsp & \longrightarrow &\End(\RR^n) \\{w} & \longrightarrow & {\Ad_w : v \longrightarrow \Ad_w(v)= w\cdot v\cdot w^{-1}.}
\end{eqnarray*}
Here $\Clrsp$ denotes the group of units of $\Clrs$. Since $\SSS^1  \cap \Spinrs =\{\pm 1\}$, we define the complex spin group by 
$$\Spinrsc = \Spinrs \times_{\ZZ_2} \SSS^1.$$
The complex spin group is the double cover of $\SOrs\times\SSS^1$, this yields to the exact sequence
$$
1 \longrightarrow \ZZ_2  \longrightarrow \Spinrsc
\stackrel{\xi^c}{\longrightarrow}  \SOrs\times\SSS^1  \longrightarrow 1,
$$
where $\xi^c =(\xi, \id^2)$. When $n=2m$ is even, $\CClrs$ has a unique irreducible complex 
representation $\chi_{2m}$ of complex dimension  $2^{m}$, $\chi_{2m} : \CClrs \longrightarrow \End(\Srs).$ If $n=2m+1$ is odd, $\CClrs$ has two inequivalent irreducible representations
both of complex dimension $2^{m},$ $\chi_{2m+1}^j : \CClrs \longrightarrow \End(\Srs^j), \ \text{for}\ j=0 \ \text{or}\ 1,$
where $\Srs^j = \{ \sigma \in \Sigma_{r,s},\ \ \chi_{2m+1}^j(\omega_{r,s})\sigma = (-1)^j\sigma\}$ 
and $\omega_{r,s}$ is the complex volume element
$$\omega_{r,s}= 
 \left\{
\begin{array}{l}
i^{m-s}\ \ \  \ \ e_1\cdot...\cdot e_n\ \ \ \ \ \ \ \ \text{if $n= 2m$,} \\
i^{m-1+s}\ \ e_1\cdot...\cdot e_n\ \ \ \ \ \ \ \text{if $n=2m+1$.}
\end{array}
\right.
$$
We define the complex spinorial representation $\rho_{n}$ by the restriction of an irreducible representation of $\CClrs$ to $\Spinrsc$:
 $$\rho_n :=\left\{
\begin{array}{l}
{\chi_{2m}}_{\arrowvert_{\Spinrsc}}\ \ \ \ \ \ \ \ \ \text{if $n =2m,$}\\ \\
{\chi_{2m+1}^0}_{\arrowvert_{\Spinrsc}}\ \ \ \ \ \ \text{if $n =2m+1,$}
\end{array}
\right.
.$$
When $n=2m$ is even, $\rho_n$ decomposes into two inequivalent irreductible representations $\rho_n^+$ and $\rho_n^-$, i.e., $\rho_n = \rho_n^+ + \rho_n^-: \Spinrsc \to \Aut(\Srs).$ The space $\Srs$ decomposes into $\Srs = \Sprs \oplus \Smrs$, where $\omega_{r,s}$ acts on $\Sprs$ as the identity and minus the identity on $\Smrs$. If $n =r+s$ is odd and when restricted to $\Spinrsc$, the two representations 
${\chi_{2m+1}^0}_{\arrowvert_{\Spinrsc}}$ and ${\chi_{2m+1}^1}_{\arrowvert_{\Spinrsc}}$ are equivalent and we
simply choose $\Srs := \Sors$. The complex spinor bundle $\Srs$ carries a Hemitian symmetric bilinear $\Spinrsc$-invariant form  $\<\cdot,
\cdot\>$, such that $$\<v\cdot\sigma_1,\sigma_2\>
= (-1)^{s+1}\<\sigma_1,v\cdot\sigma_2\> \ \text{for all}\ \  \sigma_1, \sigma_2 \in \Srs \ \text{and}\ \ v\in \RR^n.$$
Now, we give the following isomorphism $\alpha$, which is of particular importance for the identification of the $\Spinc$ bundles in the context of immersions of hypersurfaces:
\begin{eqnarray}
 \alpha: \CClrs &\to& \CClorls\nonumber \\
e_j &\to&  \nu\cdot e_j,
\label{nancy}
\end{eqnarray}
where we look at an embedding of $\RR^n$ onto $\RR^{n+1}$ such that $(\RR^n)^\perp$ is spacelike and spanned by a spacelike unit vector $\nu$.\\ \\
Let $N^n$ be an oriented semi-Riemannian manifold of signature $(r,s)$ and let $\PSON$ be the $\SOrs$-principal bundle 
of positively space and time oriented orthonormal tangent frames. A complex $\Spinc$ structure on $N$ is a $ \Spinrsc$-principal bundle $\PSpincN$ over $N$, 
an  $\SSS^1$-principal bundle $P_{\SSS^1}N$ over $N$ together with a twofold covering map $\Theta : \PSpincN \longrightarrow \PSON\times_N P_{\SSS^1}N $ such that 
$$\Theta(ua) = \Theta (u)\xi^c (a),$$
for every $u \in \PSpincN$ and $a \in \Spinrsc$, i.e., $N$ has a $\Spinc$ structure if and only if there exists an $\SSS^1$-principal bundle $P_{\SSS^1}N$ over $N$ such that the transition functions $ g_{\alpha \beta} \times l_{\alpha \beta}: U_\alpha \cap U_\beta \longrightarrow \SOrs \times \SSS^1$ of the $\SOrs \times \SSS^1$-principal bundle $\PSON\times_{N} P_{\SSS^1}N$ admit lifts to $\Spinrsc$ denoted by $\widetilde g_{\alpha \beta} \times\widetilde l_{\alpha \beta}: U_\alpha \cap U_\beta \longrightarrow  \Spinrsc,$ such that $\xi^c\circ (\widetilde g_{\alpha \beta} \times\widetilde l_{\alpha \beta}) = g_{\alpha \beta} \times l_{\alpha \beta}$. This, anyhow, is equivalent to the second Stiefel-Whitney class $w_2(N)$ being equal, modulo 2, to the first chern class $c_1 (L^N)$ of the complex line bundle $L^N$. It is the complex line bundle associated with the $\SSS^1$-principal fibre bundle via the standard representation of the unit circle.\\ \\
Let $\ScN := \PSpincN \times_{\rho_n} \Srs$ be the spinor bundle associated with the spinor representation. A section of $\ScN$ will be called a spinor field. Using the cocycle condition of the transition functions of the two principal fibre bundles $\PSpincN$ and $\PSON\times_{N} P_{\SSS^1}N$, we can prove that $$\ScN = \Sigma^{'} N \otimes (L^N)^{\frac 12},$$ where $\Sigma^{'} N$ is the locally defined spin bundle and $(L^N)^{\frac 12}$ is locally defined too but $\ScN$ is globally defined. The tangent bundle $TN = \PSON \times_{\rho_0} \RR^n$ where $\rho_0$ stands for the standard matrix representation of $\SOrs$ on $\RR^n$, can be seen as the associated vector bundle $TN \simeq \PSpincN \times_{pr_1 \circ \xi^c \circ \rho_0} \RR^n$ where $pr_1$ is the first projection. One defines the Clifford multiplication at every point $p\in N$:
\begin{eqnarray*}
T_p N  \otimes \SpcN  &\longrightarrow& \SpcN \\
\ [b,v]\otimes [b,\sigma] &\longrightarrow& [b,v]\cdot[b,\sigma]:= [b,v\cdot\sigma =\chi_n(v)\sigma],
\end{eqnarray*}
where $b\in \PSpincN$, $v\in\RR^n$, $\sigma\in\Srs$ and $\chi_n = \chi_{2m}$ if $n$ is even and $\chi_n = \chi_{2m+1}^0$ if $n$ is odd. The Clifford multiplication can be extended to differential forms. Clifford multiplication inherits the relations of the Clifford algebra, i.e., for $X,Y \in T_pN$ and $\phi\in\SpcN$ we have $X\cdot Y \cdot \phi + Y \cdot X \cdot \phi =-2 \<X,Y\> \phi.$ In even dimensions the spinor bundle splits into $\ScN = \Sigma^+ N \oplus \Sigma^- N,$ where $\Sigma^\pm N= \PSpincN \times_{\rho_n^\pm} \Sigma^\pm_{r,s}$.
Clifford multiplication by a non-vanishing tangent vector interchanges $\Sigma^+N$
and $\Sigma^- N$. The $\Spinrsc$-invariant nondegenerate symmetric sesquilinear form on $\Srs$
and $\Sigma^\pm_{r,s}$ induces inner products
on $\ScN$ and $\Sigma^\pm N$ which we again denote by $\<\cdot,\cdot\>$ and it satisfies
$$\<X\cdot\psi, \phi\> = (-1)^{s+1} \<\psi, X\cdot\phi\>,$$
for every $X\in \Gamma(TN)$ and $\psi, \phi \in \Gamma(\ScN).$ Additionally, given a connection 1-form $A^N$ on $P_{\SSS^1}N$, $A^N: T(P_{\SSS^1}N)\longrightarrow i\RR$ and the connection 1-form $\oN$ on $\PSON$ for the Levi-Civita connection $\nN$, 
we can define the connection  $$\oN\times A^N: T(\PSON\times_{N} P_{\SSS^1}N) \longrightarrow \mathfrak{so}_n\oplus i\RR =\mathfrak{spin}^{\CC}_n$$ on the principal fibre bundle $\PSON\times_{N} P_{\SSS^1}N$ and hence a covariant derivative $\ncSN$ on $\ScN$  \cite{6} given locally by 
\begin{eqnarray}
\ncSN_{e_k}\phi &= &
\Big[\widetilde{b \times s}, {e_k}(\sigma) + 
\frac{1}{4}\sum_{j=1}^n \eps_j\, e_j \cdot \nN_{e_k} e_j \cdot \sigma+ \frac 12 A^N(s_*(e_k))\sigma\Big]\nonumber\\
&=& e_k (\phi) + \frac{1}{4}\sum_{j=1}^n \eps_j\, e_j \cdot \nN_{e_k} e_j \cdot\phi + \frac 12 A^N(s_*(e_k)) \phi,
\label{zshg}
\end{eqnarray}
where $\phi = [\widetilde{b \times s},\sigma]$ is a locally defined spinor field, $ b = (e_1,\ldots,e_n)$ is a local space and time oriented orthonormal tangent frame, $s: U\longrightarrow 
P_{\SSS^1}N$ is a local section of $P_{\SSS^1}N$ and $\widetilde{b \times s}$ is the lift of the local section $b \times s: U \rightarrow \PSON\times_{N} P_{\SSS^1}N$ to the 2-fold covering $\Theta : \PSpincN \longrightarrow \PSON\times_N P_{\SSS^1}N $. The curvature of $A^N$ is an imaginary valued 2-form denoted by $F_{A^N}= dA^N$, i.e., $F_{A^N} = i\Omega^N$, where $\Omega^N$ is a real valued 2-form on $P_{\SSS^1}N$. We know
that $\Omega^N$ can be viewed as a real valued 2-form on $N$ $\cite{6}$. In this case $i\Omega^N$ is the curvature form of the associated line bundle $L^N$. 
The curvature tensor $\RSNc$ of $\ncSN$ is given by
\begin{eqnarray}
\RSNc(X,Y)\phi = 
\frac14 \sum_{j, k=1}^n \eps_j\eps_k \< \RN(X,Y)e_j,e_k\> e_j \cdot e_k \cdot\phi+ \frac i2\Omega^N(X,Y)\phi, 
\label{alphaa}
\end{eqnarray}
where $\RN$ is the curvature tensor of the Levi-Civita connection $\nN$. In the $\Spinc$ case, the Ricci identity translates, for every $X\in \Gamma(TN)$, to 
\begin{equation}
\sum_{k=1}^n \eps_k\, e_k \cdot \RSNc(e_k,X)\phi =
\half \RicN(X) \cdot \phi -\frac i2 (X\lrcorner\Omega^N)\cdot\phi.
\label{ricci}
\end{equation}
Here $\RicN$ denotes the Ricci curvature considered as a field of endomorphism
 on $TN$. The Ricci curvature considered as a symmetric bilinear form will be
written $\ricN(Y,Z) = \<\RicN(Y),Z\>$.
The  Dirac operator maps spinor fields to spinor fields and is locally defined  by
\begin{eqnarray*}
D^N\phi = i^s\sum_{j=1}^n \eps_j e_j \cdot \ncSN_{e_j} \phi, 
\end{eqnarray*}
for every spinor field $\phi$. The Dirac operator is an elliptic operator, formally selfadjoint, i.e., if $\psi$ or $\phi$ has compact support, then $(D^N\phi,\psi) = (\phi,D^N\psi)$, where $(\phi,\psi) = \int_N \<\phi,\psi\> dv_g$.
\section{Semi-Riemannian Spin$^c$ hypersurfaces and the Guass formula}
\label{hyper}
In this section, we study $\Spinc$ structures of hypersurfaces, such as the restriction of a $\Spinc$ bundle of an ambient semi-Riemannian manifold and the complex spinorial Gauss formula. \\ \\
Let $\cZ$ be an oriented ($n+1$)-dimensional semi-Riemannian $\Spinc$ manifold and $M \subset \cZ$ a semi-Riemannian hypersurface with trivial 
spacelike normal bundle. This means that there is a vector field $\nu$ on $\cZ$ along $M$ satisfying
$\<\nu,\nu\> = +1$ and $\<\nu,TM\>=0$. Hence if the signature of $M$ is $(r,s)$, then the signature of $\cZ$
is $(r+1,s)$.
\begin{prop}
The hypersurface $M$ inherts a $\Spinc$ structure from that on $\cZ$, and we have 
$$\left\{
\begin{array}{l}
\Sigma \cZ_{|_M} \ \ \ \simeq \Sigma M\ \ \ \ \text{\ \ \ if\ $n$ is even,} \\
 \SZp_{|_M} \simeq\Sigma M  \ \text{\ \ \ \ \ if\ $n$ is odd.}
\end{array}
\right.
$$
Moreover Clifford multiplication by a vector field $X$, tangent to $M$, is given by 
\begin{eqnarray}
X\bullet\phi = (\nu\cdot X\cdot \psi)_{|_M},
\label{Clifford}
\end{eqnarray}
where $\psi \in  \Gamma(\Sigma \cZ)$ (or $\psi \in \Gamma(\SZp)$ if $n$ is odd), $\phi$ is the restriction of $\psi$ to $M$, ``$\cdot$'' is the Clifford multiplication on $\cZ$, and ``$\bullet$'' that on $M$.
\label{restriction}
\end{prop}
{\bf Proof:} The bundle of space and time oriented orthonormal frames  of $M$ can be embedded into the bundle of space and time oriented orthonormal frames of $\cZ$ restricted to $M$, by
\begin{eqnarray}
 \Phi :\ \ \ \ \ \ \  \PSOM  & \longrightarrow &\PSOZ_{|_M}\\
(e_1, \cdots,e_n) & \longrightarrow & (\nu, e_1,\cdots, e_n).\nonumber
\label{phi}
\end{eqnarray}
The isomorphism $\alpha$, defined in (\ref{nancy}) yields the following commutative diagram: 
$$
\xymatrix{
\Spinrsc \ar[d]^{\xi^c}  &\hookrightarrow &\Spinrscc \ar[d]_{\xi^c} 
 \\ 
\SOrs \times \SSS^1   &\hookrightarrow & \SOrss \times\SSS^1
}
$$
where the inclusion of $\SOrs$ in $\SOrss$ is that which fixes the first basis vector under the action of $\SOrss$ on $\RR^{n+1}$. This allows to pull back via $\Phi$ the principal bundle $\PSpincZ_{|_M}$ as a $\Spinc$ structure for $M$, denoted by $\PSpincM$. Thus, we have the following commutative diagram:
$$
\xymatrix{
\PSpincM \ar[d]^{\Theta} \ar[r] &  \PSpincZ_{|_M}\ar[d]_\Theta  
& & \\
\PSOM \times_M \PSZ_{|_M} \ar[r]   & \PSOZ_{|_M}\times_M \PSZ_{|_M}
}
$$
The $\Spinrsc$-principal bundle $(\PSpincM,\pi,M)$ and the $\SSS^1$-principal bundle $(\PSM =: \PSZ|_M,\pi,M)$ define a $\Spinc$ structure on $M$.
Let $\Sigma\cZ$ be the spinor bundle on $\cZ$,
$$\Sigma\cZ = \PSpincZ\times_{\rho_{n+1}}\Sigma_{r+1,s},$$
where $\rho_{n+1}$ stands for the spinorial representation of $\Spinrscc$.
Moreover, for any spinor $\psi=[\widetilde {b\times s}, \sigma]\in \Sigma\cZ$ we can always assume that $pr_1\circ\Theta
(\widetilde {b\times s})= b$ is a local section of $\PSOZ$ with $\nu$ for first basis vector where $pr_1$ is the projection into $\PSOZ$ . Then we have 
$$\psi_{|_M} = [\widetilde{b\times s} _{|_{U\cap M}}, \sigma_{|_{U\cap M}}],$$
where the equivalence class is reduced to elements of $\Spinrsc$.
It follows that one can realise the restriction to $M$ of the spinor bundle $\Sigma\cZ$ as 
$$\Sigma\cZ_{|_M} = \PSpincM\times_{\rho_{n+1}\circ\alpha}\Sigma_{r+1,s}.$$
If $n =2m$ is even, it is easy to check that $\chi_{2m+1}^0\circ\alpha = {\chi_{2m+1}^0}_{\arrowvert_{\CClorls}}.$ Hence $\chi_{2m+1}^0\circ\alpha$ is an irreductible representation of $\CClrs$ of dimension $2^{m}$, as ${\chi_{2m+1}^0}_{\arrowvert_{\CClorls}}$, and finally
$\chi_{2m+1}^0\circ\alpha \cong \chi_{2m}.$ We conclude that
$$ \rho_{2m+1}\circ\alpha \cong \rho_{2m},\ \ \ \text{and}\ \ \ \Sigma \cZ_{|_M}\cong \Sigma M.$$
If $n = 2m+1$ is odd, we know that $\chi_{2m+1}^0$ is the unique irreductible representation of  $\CClrs$ of dimension $2^m$ for which the action of the complex volume form is the identity. Since $n+1 = 2m+2$ is even, $\SZ$ decomposes into positive and negative parts, $\SZpm = \PSpincZ\times_{\rho_{2m+1}^\pm}\Sigma_{r+1,s}^\pm.$  It is easy to show that 
$\chi_{2m+2}\circ\alpha = {\chi_{2m+2}}_{\arrowvert_{\CClorls}},$ but $\chi_{2m+2}\circ\alpha$ can be written as the direct sum of two irreductible inequivalent representations, as 
${\chi_{2m+2}}_{\arrowvert_{\CClorls}}$. Hence, we have 
$$\chi_{2m+2}\circ\alpha=(\chi_{2m+2}\circ\alpha)^+\oplus(\chi_{2m+1}\circ\alpha)^-,$$
where  $(\chi_{2m+2}\circ\alpha)^\pm (\omega_{r,s}) = \pm \id_{\Srs}$. The representation $\chi_{2m+1}^0$ being the unique representation of $\CClrs$ of dimension $2^m$ for which the action of the volume form is the identity, we get
$(\chi_{2m+2}\circ\alpha)^+  \cong  \chi_{2m+1}^0.$ Finally,
$$\rho_{2m+2}^+\circ\alpha  \cong  \rho_{2m+1}\ \ \ \ \text{and}\ \ \ \ \SZp_{|_M} \cong \Sigma M.$$
Now, Equation (\ref{Clifford}) follows directly from the above identification.
\begin{remarks}
 $1.$ The algebraic remarks in the previous section show that if $n$ is odd we can also get $\SZm_{|_M} \simeq \Sigma M,$
where the Clifford multiplication by a vector field tangent to $M$ is given by 
$X\bullet\phi = -(\nu\cdot X\cdot\psi)_{|_M}.$\\
$2.$ The connection 1-form defined on the restricted $\SSS^1$-principal bundle $(\PSM =: \PSZ_{|_M},\pi,M)$, is given by
$$A^M = {A^\cZ}_{|_M} : T(\PSM) = T(\PSZ)_{|_M} \longrightarrow i\RR.$$
Then the curvature 2-form $i\Omega^M$ on the $\SSS^1$-principal bundle $\PSM$ is given by $i\Omega^M = {i\Omega^\cZ}_{|_M}$, which can be viewed as an imaginary 2-form on $M$ and hence as the curvature form of the line bundle $L^M$, the restriction of the line bundle $L^\cZ$ to $M$.\\
$3.$  For every $\psi \in \Gamma(\SZ)$ ($\psi \in \Gamma(\SZp)$ if $n$ is odd), the real 2-forms $\Omega^M$ and $\Omega^\cZ$ are related by the following formulas:
\begin{eqnarray}
|\Omega^\cZ|^2 = |\Omega^M|^2 + |\nu\lrcorner\Omega^\cZ|^2,
 \label{moduleomegareltion}
\end{eqnarray}
\begin{eqnarray}
(\Omega^\cZ \cdot\psi)_{|_M} = \Omega^M\bullet\phi + (\nu\lrcorner\Omega^\cZ)\bullet\phi.
\label{glucose}
\end{eqnarray}
In fact, we can write
\begin{eqnarray*}
\Omega^\cZ = \sum_{i=1}^n \Omega^\cZ (\nu, e_i)\ \nu\wedge e_i + \sum_{i<j}^n \Omega^\cZ \ (e_i, e_j)\ e_i\wedge e_j=-(\nu \lrcorner\Omega^\cZ) \wedge \nu + \Omega^M,
\end{eqnarray*}
which is (\ref{moduleomegareltion}). When restricting the Clifford multiplication of $\Omega^\cZ$ by $\psi$ to the hypersurface $M$ we obtain
\begin{eqnarray}
(\Omega^\cZ\cdot\psi)_{|_M} = \big(\nu\cdot(\nu\lrcorner\Omega^\cZ)\cdot\psi\big)_{|_M} + (\Omega^M\cdot\psi)_{|_M}
= (\nu\lrcorner\Omega^\cZ)\bullet\phi + \Omega^M\bullet\phi.
\end{eqnarray}
\end{remarks}
\begin{prop} [The spinorial Gauss formula]
We denote by $\nabla^{\Sigma \cZ}$ the spinorial Levi-Civita connection on $\Sigma \cZ$ and by $\nabla^{\Sigma M}$ that on $\Sigma M$. For all $X\in \Gamma(TM)$ and for every spinor field $\psi \in \Gamma(\Sigma \cZ)$, then
\begin{equation}
(\nabla^{\Sigma \cZ}_X\psi)_{|_M} =  \nabla^{\Sigma M}_X \phi - \half W(X)\bullet\phi,
\label{spingauss}
\end{equation}
where $W$ denotes the Weingarten map with respect to $\nu$ and $\phi= \psi_{|_M}$. Moreover, let $D^\cZ$ and $D^M$ be the Dirac operators on $\cZ$ and $M$. Denoting by the same symbol any spinor and it's restriction to $M$, we have
\begin{equation}
\nu\cdot D^\cZ\phi = \DDM\phi +\frac{i^s n}{2}H\phi - i^s \nabla^{\Sigma\cZ}_{\nu}\phi,
\label{diracgauss}
\end{equation}
where $H = \frac1n \tr(W)$ denotes the mean curvature and 
$\DDM = D^M$ if $n$ is even and $\DDM=D^M \oplus(-D^M)$ if $n$ is odd.
\end{prop}
{\bf Proof:} The Riemannian Gauss formula is given, for every vector fields $X$ and $Y$ on $M$, by 
\begin{equation}
\nZ_XY = \nM_XY + \<W(X),Y\>\nu.
\label{weingarten}
\end{equation}
Let $(e_1, e_2,...,e_n)$ a local space and time oriented orthonormal frame of $M$, such that $ b= (e_0= \nu, e_1, e_2,...,e_n)$ is that of $\cZ$. We consider $\psi$  a local section of $\Sigma\cZ$, $\psi = [\widetilde{b\times s},\sigma]$ where $s$ is a local section of $\PSZ$. Using (\ref{zshg}), (\ref{weingarten}) and the fact that $X(\psi)|_M = X(\phi)$ for $X\in \Gamma(TM)$, we compute for $j=1,...,n$
\begin{eqnarray*}
\Big(\nabla_{e_j}^{\Sigma\cZ}\psi\Big)_{|_M} &=&
e_j(\phi) + \frac 14 \sum_{k=0}^n\eps_k (e_k\cdot\nabla_{e_j}^{\cZ}e_k \cdot\psi)_{|_M} + \frac 12  A^Z (s_*(e_j))\phi \\&=&
e_j(\phi) + \frac 14 \sum_{k=1}^n\eps_k (e_k\cdot\nabla_{e_j}^{\cZ}e_k \cdot\psi)_{|_M}+\frac 14 (\nu\cdot\nabla_{e_j}^{\cZ}\nu\cdot\psi)_{|_M}+\frac12  A^M (s_*(e_j))\phi \\
&=&\nabla_{e_j}^{\Sigma M}\phi + \frac 14 \sum_{k=1}^n \eps_k <W(e_j), e_k>(e_k\cdot\nu\cdot\psi)_{|_M}- \frac 14 (\nu\cdot W(e_j)\cdot\psi)_{|_M}\\
&=& \nabla_{e_j}^{\Sigma M}\phi - \half (\nu \cdot W(e_j) \cdot \psi)_{|_M}\\ &=&  \nabla^{\Sigma M}_{e_j} \phi - \half W(e_j)\bullet\phi.
\end{eqnarray*}
Moreover $(D^\cZ\psi)_{|_M} = i^{s}\sum_{j=1}^n \eps_j (e_j \cdot \nabla^{\Sigma\cZ}_{e_j} \psi)_{|_M}
+i^{s} (\nu \cdot \nabla^{\Sigma\cZ}_{\nu} \psi)_{|_M}$, and by (\ref{spingauss}),
\begin{eqnarray*}
i^{s}\sum_{j=1}^n \eps_j (e_j \cdot \nabla^{\Sigma \cZ}_{e_j} \psi)_{|_M}
&=&
i^{s}\sum_{j=1}^n \eps_j\, (e_j \cdot \nabla^{\Sigma M}_{e_j} \phi)
-i^{s}\half \sum_{j=1}^n \eps_j\, (e_j \cdot \nu \cdot W(e_j) \cdot \psi)_{|_M} \\
&=&
- i^{s}\nu\cdot \sum_{j=1}^n \eps_j\, \nu\cdot e_j \cdot \nabla^{\Sigma M}_{e_j} \phi  
+i^{s}\half \sum_{j=1}^n \eps_j\, (\nu \cdot e_j \cdot  W(e_j) \cdot \psi)_{|_M}\\
&=&
-\nu\cdot \DDM \phi-  \frac{i^{s}}{2}\tr(W) (\nu \cdot\psi)_{|_M}.
\end{eqnarray*}
\begin{prop}
Let $\cZ$ be an $(n+1)$-dimensional semi-Riemannian $\Spinc$ manifold.
Assume that $\cZ$ carries a semi-Riemannian foliation by hypersurfaces with trivial
spacelike normal bundle, i.e., the leaves $M$ are semi-Riemannian 
hypersurfaces and there exists a vector field $\nu$ on $\cZ$ perpendicular
to the leaves such that $\<\nu,\nu\>=1$ and $\nZ_\nu\nu = 0$. Then the commutator of the leafwise Dirac operator and the normal
derivative is given by
$$
i^{-s}[\nabla^{\Sigma\cZ}_\nu,\DDM]\,\phi = 
\DW\phi - \frac{n}{2}\, \nu\cdot\gradM(H)\cdot\phi 
+ \half\,\nu\cdot\divM(W)\cdot\phi +\frac i2 \nu\cdot(\nu\lrcorner\Omega^\cZ)\cdot\phi  .
$$
Here $\gradM$ denotes the leafwise gradient, $\divM(W) = \sum_{i=1}^n
\eps_i\, (\nM_{e_i}W)(e_i)$ denotes the leafwise divergence of the
endomorphism field $W$ and  $\DW\phi =
\sum_{i=1}^n\eps_i\, \nu\cdot e_i\cdot\nabla^{\Sigma M}_{W(e_i)}\phi$.
\label{kommutator}
\end{prop}
{\bf Proof:} We choose a local oriented orthonormal tangent frame $(e_1,\ldots,e_n)$ 
for the leaves and we may assume for simplicity that $\nZ_\nu e_j = 0$. Now, we compute
\begin{eqnarray}
i^{-s}[\nabla^{\Sigma \cZ}_\nu,\DDM]\,\phi &=&
\summe \eps_j \left( \nabla^{\Sigma \cZ}_\nu(\nu \cdot e_j\cdot \nabla^{\Sigma M}_{e_j}\phi) 
-  \nu \cdot e_j\cdot \nabla^{\Sigma M}_{e_j} \nabla^{\Sigma \cZ}_\nu\phi  \right)
\nonumber\\
&=&
\summe \eps_j\,\nu \cdot e_j\cdot  \left( \nabla^{ \Sigma\cZ}_\nu\nabla^{\Sigma M}_{e_j}\phi 
- \nabla^{\Sigma M}_{e_j} \nabla^{\Sigma \cZ}_\nu\phi  \right)
\nonumber\\
&\stackrel{(\ref{spingauss})}{=}&
\summe \eps_j\,\nu \cdot e_j\cdot  \Big[ \nabla^{\Sigma \cZ}_\nu(\nabla^{\Sigma \cZ}_{e_j}+
\half \nu\cdot W(e_j))
\nonumber\\
&&
-  (\nabla^{\Sigma \cZ}_{e_j}+\half \nu\cdot W(e_j)) \nabla^{\Sigma \cZ}_\nu \Big]\phi 
\nonumber\\
&=&
\summe \eps_j\,\nu \cdot e_j\cdot \Big(\RSZc(\nu,e_j)+\nabla^{\Sigma\cZ}_{[\nu,e_j]}
+ \half \nu\cdot (\nZ_\nu W)(e_j)   \Big)\phi
\nonumber\\
&\stackrel{(\ref{ricci})}{=}&
-\half\nu\cdot\RicZ(\nu)\cdot\phi +\frac i2 \nu\cdot(\nu\lrcorner\Omega^\cZ)\cdot\phi\nonumber
\\ &\ & +\summe \eps_j\,\nu \cdot e_j\cdot \Big(\nabla^{\Sigma \cZ}_{W(e_j)}
+ \half \nu\cdot (\nZ_\nu W)(e_j)\Big)\phi
\nonumber\\
&\stackrel{(\ref{spingauss})}{=}&
-\half\nu\cdot\RicZ(\nu)\cdot\phi+\frac i2 \nu\cdot(\nu\lrcorner\Omega^\cZ)\cdot\phi 
\nonumber\\
&&
+ \summe \eps_j\,\nu \cdot e_j\cdot \Big(\nabla^{\Sigma M}_{W(e_j)} -\half\nu\cdot W^2(e_j)
+ \half \nu\cdot (\nZ_\nu W)(e_j)   \Big)\phi
\nonumber\\
&=&
-\half\nu\cdot\RicZ(\nu)\cdot\phi +\frac i2 \nu\cdot(\nu\lrcorner\Omega^\cZ)\cdot\phi+\DW\phi
\nonumber\\
&&
+ \half \summe \eps_j\, e_j\cdot \Big(-W^2(e_j)
+ (\nZ_\nu W)(e_j)   \Big)\phi\nonumber .
\end{eqnarray}
The Riccati equation for the Weingarten map $(\nZ_\nu W)(X) = 
\RZ(X,\nu)\nu + W^2(X)$ yields
\begin{eqnarray}
i^{-s}[\nabla^{\Sigma \cZ}_\nu,\DDM]\,\phi &=&
-\half\nu\cdot\RicZ(\nu)\cdot\phi +\frac i2 \nu\cdot(\nu\lrcorner\Omega^\cZ)\cdot\phi+\DW\phi
\nonumber\\ &&+ \half \summe \eps_j\, e_j\cdot (\RZ(e_j,\nu)\nu)\cdot\phi
\nonumber\\
&=&
-\half\nu\cdot\RicZ(\nu)\cdot\phi +\frac i2 \nu\cdot(\nu\lrcorner\Omega^\cZ)\cdot\phi+\DW\phi
+ \half \ricZ(\nu,\nu)\phi
\nonumber\\
&=&
\DW\phi -\half \summe \eps_i\, \ricZ(\nu,e_j)\, \nu\cdot e_j\cdot\phi+\frac i2 \nu\cdot(\nu\lrcorner\Omega^\cZ)\cdot\phi.
\label{ausdruckmitricci}
\end{eqnarray}
The Codazzi-Mainardi equation for $X,Y,V
\in TM$ is given by $\<\RZ(X,Y)V,\nu\> = \<(\nM_X W)(Y),V\> - \<(\nM_Y W)(X),V\> .$ Thus,
\begin{eqnarray*}
\ricZ(\nu,X) &=&
\summe \eps_j \<\RZ(X,e_j)e_j,\nu\> \\ &=& 
\summe \eps_j \left( \<(\nM_X W)(e_j),e_j\> - \<(\nM_{e_j} W)(X),e_j\>\right)
\nonumber\\
&=&
\tr(\nM_X W) - \<\divM(W),X\> .
\end{eqnarray*}
Plugging this into (\ref{ausdruckmitricci}) we get
\begin{eqnarray*}
i^{-s}[\nabla^{\Sigma \cZ}_\nu,\DDM]\,\phi 
&=&
\DW\phi -\half \summe \eps_j\, \left(\tr(\nM_{e_j} W) 
- \<\divM(W),e_j\>\right)\, \nu\cdot e_j\cdot\phi \\ &&+\frac i2 \nu\cdot(\nu\lrcorner\Omega^\cZ)\cdot\phi.\\
&=&
\DW\phi -\half \summe \eps_j\, e_j (\tr(W)) \nu\cdot e_j\cdot\phi
+\half \nu\cdot\divM(W) \cdot\phi \\ &&+ \frac i2 \nu\cdot(\nu\lrcorner\Omega^\cZ)\cdot\phi.\\
&=&
\DW\phi -\frac{n}{2}\nu\cdot \gradM(H)\cdot\phi 
+\half \nu\cdot\divM(W) \cdot\phi+\frac i2 \nu\cdot(\nu\lrcorner\Omega^\cZ)\cdot\phi. 
\end{eqnarray*}
\section{The generalized cylinder on semi-Riemannian Spin$^c$ manifolds}
\label{gc}
Let $M$ be an $n$-dimensional smooth manifold and $g_t$ a smooth 1-parameter family of semi-Riemannian metrics on $M$,
$t\in I$ where $I \subset \RR$ is an interval.
We define the  generalized cylinder by $$\cZ := I \times M,$$ with semi-Riemannian metric $g_\cZ := \<\cdot,\cdot\> =  dt^2 + g_t$.
The generalized cylinder is an $(n+1)$-dimensional semi-Riemannian
manifold of signature $(r+1,s)$ if the signature of $g_t$ is $(r,s)$.
\begin{prop}
There is a 1-1-correspondence between the $\Spinc$ structures on $M$ and that on $\cZ$.
\label{omegazero}
\end{prop}
{\bf Proof: }
As explained in Section \ref{hyper}, $\Spinc$ structures on $\cZ$ can be 
restricted to $\Spinc$ structures on $M$.
Conversely, given a $\Spinc$ structure on $M$ it can be pulled back
to $I \times M$ via the projection $pr_2 : I \times M \longrightarrow M$ yields a $\Spinc$ structure on $\cZ$. In fact, the
pull back of the $\Spinrsc$-principal bundle  $\PSpincM$ on $M$ gives rise to a $\Spinrsc$-principal bundle on $\cZ$ denoted by  $\PSpincZ$
$$
\xymatrix{
\PSpincZ \ar[d]^{\pi} \ar[r] & \PSpincM \ar[d]_\pi 
& & \\
\cZ = I \times M \ar[r]  & M  
}
$$
Enlarging the structure group via the embedding  $\Spinrsc 
\hookrightarrow
\Spinrscc $, which covers the standard embedding
\begin{eqnarray*}
\SOrs\times \SSS^1 &\hookrightarrow&
\SOrss\times \SSS^1\\
(a,z)&\mapsto& \Big(\Big(
\begin{tabular}{ll}
1 & 0 \\ 
0 & a
\end{tabular}
\Big), z\Big),
\end{eqnarray*}
gives a $\Spinrscc$-principal fibre bundle on $\cZ$, denoted also by  $\PSpincZ$. The pull back of the line bundle $L^M$ on $M$ defining the $\Spinc$ structure on $M$, gives a line bundle $L^\cZ$ on $\cZ$
such that the following diagram commutes
$$
\xymatrix{
L^\cZ=pr_2^*(L^M)\ar[d]^{\pi} \ar[r] & L^M \ar[d]_\pi 
& &\\
\cZ = I \times M \ar[r]  & M  
}
.$$
The line bundle $L^\cZ$ on $\cZ$ and the $\Spinrscc$-principal fibre bundle $\PSpincZ$ on $\cZ$  yields the $\Spinc$ structure on $\cZ$ which restricts to the given $\Spinc$
structure on $M$.
\begin{remark}
If $M$ is a $\Spinc$ Riemannian manifold and if we denote by $i\Omega^M$ the imaginairy valued curvature on the line bundle $L^M$, we know that there  exists a unique curvature 2-form, denoted by $i\Omega^\cZ$, on the line bundle $L^\cZ = pr_2 ^* (L^M)$, defined by 
$i\Omega^\cZ = pr_2^* (i\Omega^M).$ Thus we have $$\Omega^\cZ(X, Y) = \Omega^M (X, Y)\ \ \text{and}\ \ \Omega^\cZ(\nu, Y) =0\ \ \text{for any}\ \ X, Y \in \Gamma(TM).$$
\label{zero}
\end{remark}
\vspace{-1.5cm}
\begin{prop} $\cite{3}$
 On a generalized cylinder $\cZ = I \times M$ with semi-Riemannian metric $g^\cZ = 
\<\cdot,\cdot\> = dt^2 + g_t$ we define, in every $p\in M$ and $X,Y \in T_pM$, the first and second derivatives of $g_t$ by 
$$\gdt(X,Y) := \frac{d}{dt}(g_t(X,Y))\ \ \ \ \text{and}\ \ \ \ 
\gddt(X,Y) := \frac{d^2}{dt^2}(g_t(X,Y)).$$ Hence the following formulas hold:
\begin{eqnarray}
\< W(X),Y\> &=& -\half \gdt(X,Y) ,
\label{zylweingarten}\\
\<R^\cZ(U,V)X,Y\> &=& \<R^{M_t}(U,V)X,Y\>
\label{zylgauss}\\ 
&& + \frac{1}{4} \Big(
\gdt(U,X) \gdt(V,Y) - \gdt(U,Y) \gdt(V,X)\Big) ,
\nonumber \\
\<R^\cZ(X,Y)U,\nu\> &=& \half \Big((\nMt_Y\gdt)(X,U) - (\nMt_X\gdt)(Y,U)\Big),
\label{zylcodazzi}\\
\<R^\cZ(X,\nu)\nu,Y\> &=& -\half\Big( \gddt(X,Y) + \gdt(W(X),Y)  \Big),
\label{zylriccati} 
\end{eqnarray}
where $X,Y,U,V \in T_pM$, $p\in M$.
\label{zylformeln}
\end{prop}
\section{The variation formula for the Dirac operator on Spin$^c$ manifolds}
\label{vf}
First we give some facts about parallel transport on $\Spinc$ manifolds along a curve $c$. We consider a Riemannian $\Spinc$ manifold $N$, we know that there exists a unique correspondence which associates to a spinor field $\psi (t)= \psi(c(t))$  along a curve $c: I\longrightarrow N$ another spinor field $\frac{D}{dt}\psi$ along $c$, called the covariant derivative of $\psi$ along $c$, such that
$$\frac{D}{dt}(\psi+\phi)=\frac{D}{dt}\psi +\frac{D}{dt}\phi,\ \ \text{for any}\ \psi \ \text{and}\ \ \phi\ \text{along the curve $c$},$$
$$\frac{D}{dt}(f\psi)= f\frac{D}{dt}\psi + (\frac{d}{dt}f)\ \psi,\text{\ where  $f$ is a differentiable function on $I$},$$
$$\nabla^{\Sigma N}_{\dot{c}(t)} \psi = \frac{D}{dt}\phi,\text{\ where $\phi(t)=\psi(c(t))$}.$$ 
A spinor field $\psi$ along a curve $c$ is called parallel when $\frac{D}{dt}\psi(t)=0$ for all $t\in I$. Now, if $\psi_0$ is a spinor at the point $c(t_0)$, $t_0 \in I, (\psi_0\in \Sigma_{c(t_0)} N)$ then there exists a unique parallel spinor $\phi$ along $c$, such that $\psi_0= \phi(t_0)$. The linear isometry $\tau_{t_0}^{t_1}$ defined by
\begin{eqnarray*}
 \tau_{t_0}^{t_1}:\Sigma_{c(t_0)} N &\longrightarrow& \Sigma_{c(t_1)} N \\
\psi_0 &\longrightarrow&  \phi(t_1),
\end{eqnarray*}
is called the parallel transport along the curve $c$ from $c(t_0)$ to $c(t_1)$. The basic property of the parallel transport on a $\Spinc$ manifold is the following: Let  $\psi$ be a spinor field on a Riemannian $\Spinc$ manifold $N$ , $X\in \Gamma(TN)$, $p \in N$ and $c:I \longrightarrow N$ an integral curve through $p$, i.e., $c(t_0)= p$ and $\frac{d}{dt}c(t)=X(c(t))$, we have
\begin{eqnarray}
(\nabla^{\Sigma N}_{X}\psi)_p = \frac {d}{dt}\Big(\tau_{t}^{t_0}(\psi(t))\Big)|_{t=t_0}.
\label{transport}
\end{eqnarray}
\\
Now, we consider $g_t$ a smooth 1-parameter family of semi-Riemannian metrics on a 
$\Spinc$ manifold $M$ and the generalized cylinder  $\cZ = I \times M$ with semi-Riemannian metric $g^\cZ = 
\<\cdot,\cdot\> = dt^2 + g_t$. For $t\in I$ we denote by $M_t$ the manifold $(M, g_t)$. Let us write  ``$\cdot$'' for the Clifford multiplication on $\cZ$ and ``$\cdott$'' for that on $\Mt$. Recall from Section \ref{gc} that $\Spinc$ structures on $M$ and $\cZ$ are in 1-1-correspondence and $\Sigma\cZ|_{\Mt} = 
\Sigma\Mt$ as hermitian vector bundles if $n=r+s$ is even and 
$\Sigma^+\cZ|_\Mt = \Sigma\Mt$ if $n$ is odd. For a given $x\in M$ and $t_0, t_1 \in I$, parallel transport $\tau_{t_0}^{t_1}$ on the generalized cylinder $\cZ$ along the curve $c: I\rightarrow I\times M, t\rightarrow (t, x)$ is given by 
$$\tau_{t_0}^{t_1}:\Sigma_{c(t_0)} \cZ \simeq \Sigma_x M_{t_0} \longrightarrow \Sigma_{c(t_1)} \cZ \simeq \Sigma_x M_{t_1}.$$
This isomorphism satisfies
$$\tau_{t_0}^{t_1} (X\bullet_{t_0}\varphi) = (\zeta_{t_0}^{t_1} X)\bullet_{t_1} (\tau_{t_0}^{t_1} \varphi),$$
$$ < \tau_{t_0}^{t_1} \psi, \tau_{t_0}^{t_1} \varphi> = <\psi, \varphi>,$$
where $\zeta_{t_0}^{t_1}: T_{(x, t_0)} \cZ \simeq T_x M_{t_0} \rightarrow T_{(x, t_1)} \cZ \simeq T_x M_{t_1}$ is the parallel transport on $\cZ$ along the same curve $c$, $X\in T_x M_{t_0}$ and $\psi, \varphi \in \Sigma_x M_{t_0}$.
\begin{thm}
\label{Diracvariation}
On a $\Spinc$ manifold $M$, let $g_t$ be a smooth 1-parameter family of semi-Riemannian metrics. Denote by $\DDMt$ the Dirac operator of $\Mt$, and $\Dgdt = \sum_{i,j=1}^n \eps_i\eps_j 
\gdt(e_i,e_j)e_i\cdott\nabla_{e_j}^{\Sigma\Mt}$.
Then for any smooth spinor field $\psi$ on $\Mto$ we have
$$
\left.\frac{d}{dt}\right|_{t=t_0} \tau_t^{t_0}\DDMt \tau_{t_0}^t \psi
=
-\half\mathfrak{D}^{\dot{g}_{t_0}}\psi + \frac{1}{4}\, \grad^{\Mto}
(\tr_{g_{t_0}}(\dot{g}_{t_0}))\bullet_{t_0}\psi 
- \frac14\,\div^{\Mto}(\dot{g}_{t_0})\bullet_{t_0}\psi .
$$
\end{thm}
{\bf Proof:} The vector field $\nu := \pa{}{t}$ is spacelike of unit length and
orthogonal to the hypersurfaces $M_t := \{t\}\times M$. Denote by $W_t$ the Weingarten map of $M_t$ with respect to $\nu$
and by $H_t$ the mean curvature. If $X$ is a local coordinate field on $M$, then $\<X,\nu\> = 0$ and
$[X,\nu] = 0$. Thus
\begin{eqnarray*}
0 &=& 
d_\nu \<X,\nu\> = \<\nZ_\nu X,\nu\> + \<X,\nZ_\nu\nu\>
= 
\<\nZ_X \nu,\nu\> + \<X,\nZ_\nu\nu\>\\ 
&=&
-\<W_t(X),\nu\> + \<X,\nZ_\nu\nu\> = \<X,\nZ_\nu\nu\>
\end{eqnarray*}
and differentiating $\<\nu,\nu\>=1$ yields $\<\nu,\nZ_\nu\nu\>=0$. Hence $\nZ_\nu\nu = 0$, i.e.,  for $x\in M$ the curves $t\mapsto (t,x)$ are geodesics parametrized
by arclength.
So the assumptions of Proposition \ref{kommutator} are satisfied for the
foliation $(M_t)_{t\in I}$. 
By Remark \ref{zero}, the commutator formula of Proposition \ref{kommutator}
gives for a section $\phi$ of $\Sigma M_t$, (or $\Sigma^+ M_t$ if $n$ is odd)
\begin{eqnarray}
i^{-s}[\nabla^{\Sigma\cZ}_\nu,\DDMt]\,\phi = 
\DWt\phi - \frac{n}{2}\, \gradMt(H_t)\cdott\phi 
&+& \half\,\divMt(W_t)\cdott\phi.
\label{kommneu}
\end{eqnarray}
From Proposition \ref{zylformeln} we deduce 
$$\divMt(W_t) = -\half \divMt(\gdt),\ \ \ \ \ \ H_t = -\frac{1}{2n} \tr_{g_t}(\gdt)\ \ \ \ \ \text{and}\ \ \ \ \ \DWt = -\half \Dgdt.$$
Thus (\ref{kommneu}) can be rewritten as 
\begin{eqnarray}
i^{-s}[\nabla^{\Sigma\cZ}_\nu,\DDMt]\,\phi = 
-\half\Dgdt\phi + \frac{1}{4}\, \gradMt(\tr_{g_t}(\gdt))\cdott\phi -\frac14\,\divMt(\gdt)\cdott\phi .
\label{kommneuneu}
\end{eqnarray}
Now if $\phi$ is parallel along the curves $t \mapsto (t,x)$, i.e., it
is of the form $\phi(t,x) = \tau_{t_0}^t \psi(t_0,x)$ for some spinor field
$\psi$ on $\Mto$, then using (\ref{transport}) at $t=t_0$, the left hand side of (\ref{kommneuneu}) could be written as
\begin{eqnarray}
i^{-s} [\nabla^{\Sigma\cZ}_\nu,\DDMt]\,\phi &=& i^{-s} \nabla^{\Sigma\cZ}_\nu\DDMt\,\phi 
= i^{-s}\left.\frac{d}{dt}\right|_{t=t_0} \tau_t^{t_0}\DDMt\,\phi \nonumber \\ &=& i^{-s}\left.\frac{d}{dt}\right|_{t=t_0} \tau_t^{t_0}\DDMt \tau_{t_0}^t \psi,
\label{youpi}
\end{eqnarray}
which gives the variation formula for the Dirac operator.
\begin{cor}
Let $(M^n,g)$ be a $\Spinc$ Riemannian manifold, if we consider the family of metrics defined by $g_t = g + tk,$
where $k$ is a symmetric $(0,2)$-tensor, we have
\begin{equation}
 \left.\frac{d}{dt}\right|_{t=0} \tau_t^{0}\DDMt \tau_{0}^t \psi
=
-\frac{1}{2}\mathfrak{D}^{k}\psi + \frac{1}{4}\, \grad^{M}
(\tr_{g}(k))\cdot\psi 
- \frac{1}{4}\,\div^{M}(k)\cdot\psi,
\label{new}
\label{fr}
\end{equation}
where ``$\cdot =\bullet_{t_0 =0} $'' is the Clifford multiplicaton on $M$.
\end{cor}
This formula has been proved in \cite{1}, Theorem $21$ for spin Riemannian manifolds and in \cite{3} for spin semi-Riemannian manifolds.
\section{Energy-Momentum tensor on Spin$^c$ manifolds}
In this section we study the Energy-Momentum tensor on $\Spinc$ Riemannian manifolds from a geometric point of vue. We begin by giving the proofs of Proposition \ref{endomo}, Theorem \ref{frkim} and Proposition \ref{morel}.\\ \\
{\bf Proof of Proposition \ref{endomo} :} Using Equation (\ref{new}) we calculate
\begin{eqnarray*}
 \left.\frac{d}{dt}\right|_{t=0} (\tau_t^0\DDMt \tau_{0}^t \psi, \psi)_{g_t} &=& 
\left.\frac{d}{dt}\right|_{t=0} (\DDMt \tau_{0}^t \psi, \tau_{0}^t \psi)_{g_t}\\
&=& -\frac 12  (\mathfrak D^k \psi,\psi)_{g}\\
&=& - \frac 12 \sum_{i,j}
k(e_i,e_j)(e_i \cdot\nabla^{\Sigma M}_{e_j}\psi,\psi)\\
&=& -\frac 12\int_M < k, T_\psi >dv_g.
\end{eqnarray*} 
{\bf Proof of Theorem \ref{frkim}} : The Proof of this Theorem will be omitted since it is similar to the one given by  
Friedrich and Kim in \cite{4} for spin manifolds.\\ \\
{\bf Proof of Proposition \ref{morel}} : Let $\psi$ be any parallel spinor field on $\cZ$. Then Equation (\ref{spingauss}) yields
\begin{eqnarray}
\nabla^{\Sigma M}_X\phi = \frac 12 W(X)\bullet\phi.
\end{eqnarray}
Let $(e_1,..., e_n)$ be a positively oriented local orthonormal basis of $TM$. For $j=1,...,n$ we have 
$$\nabla^{\Sigma M}_{e_j}\phi = \frac 12\sum_{k=1}^n  W_{jk}\ e_k\bullet\phi.$$
Taking Clifford multiplication by $e_i$ and the scalar product with $\phi$, we get
$$\Re (e_i\bullet\nabla^{\Sigma M}_{e_j}\phi, \phi) = \frac 12 \sum_{k=1}^n W_{jk} \Re (e_i\bullet e_k\bullet\phi, \phi).$$
Since $\Re(e_i\bullet e_k\bullet\phi, \phi) = -\delta_{ik} \vert\phi\vert^2$, it follows, by the symmetry of $W$
$$\Re (e_i\bullet \nabla^{\Sigma M}_{e_j}\phi + e_j\bullet\nabla^{\Sigma M}_{e_i}\phi, \phi) = -W_{ij} \vert\phi\vert^2.$$
Therefore, $2 \ell^\phi = -W$. Using Equation (\ref{diracgauss}) it is easy to see that $\phi$ is an eigenspinor associated with the eigenvalue $-\frac n2 H$ of $\widetilde D$. Since $\scalZ = \scalM + 2\ \ricZ(\nu, \nu) - n^2 H^2 + \vert W\vert^2$ we get
\begin{eqnarray*}
 \frac 14 (\scalM -c_n \vert\Omega^M\vert) + \vert T^\phi\vert^2 &=& \frac 14(\scalZ -2\ \ricZ(\nu, \nu) -c_n\vert\Omega^M\vert)+n^2 \frac{H^2}{4} \\&=& n^2 \frac{H^2}{4},
\end{eqnarray*}
hence $M$ satisfies the equality case in (\ref{rodgy}) if and only if (\ref{maga}) holds.
\begin{cor} 
Under the same conditions as Proposition \ref{morel}, if $n=2$ or $3$, the hypersurface $M$ satisfies the equality case in (\ref{rodgy}) if $\RicZ(\nu)=0$ and $\scalZ \geqslant 0$.\label{cococo}
\end{cor}
{\bf Proof:} Since $\cZ$ has a parallel spinor, we have (see \cite{6})
\begin{eqnarray}
|\RicZ (\nu)| = |\nu\lrcorner\Omega^\cZ|,
\label{111}
\end{eqnarray}
\begin{eqnarray}
 i (Y\lrcorner \Omega^\cZ)\cdot\psi = \RicZ (Y)\cdot\psi\ \ \text{for every}\ \ \ Y\in \Gamma(\Sigma\cZ).
\label{222} 
\end{eqnarray}
For $Y= e_j$ in Equation (\ref{222}) then taking Clifford multiplication by $e_j$ and summing from $j=1,..., n+1$, we get
$$i \sum_{j=1}^{n+1} e_j \cdot (e_j\lrcorner\Omega^\cZ)\cdot\psi = \sum_{j=1}^{n+1} e_j\cdot\RicZ (e_j)\cdot\psi = -\scalZ \psi.$$
But $2\ \Omega^\cZ\cdot\psi = \sum_{j=1}^{n+1} e_j \cdot (e_j\lrcorner\Omega^\cZ)\cdot\psi,$ hence we deduce that $\Omega^\cZ\cdot\psi = i\frac{\scalZ}{2}\psi$. By (\ref{111}) and (\ref{glucose}) we obtain $\Omega^M\bullet\phi = i \frac {\scalZ}{2} \phi$. Since  $n=2$ or $3$ we have $|\Omega^M| =\frac{\scalZ}{2}$ and  Equation (\ref{maga}) is satisfied.
%
\begin{cor}
 Under the same conditions as Proposition \ref{morel}, if the restriction of the complex line bundle $L^\cZ$ is flat, i.e., $L^M$ is a flat complex line  bundle ($\Omega^M = 0$), the hypersurface $M$ is a limiting manifold for (\ref{rodgy}).
\end{cor}
{\bf Proof:} Since $\Omega^M=0$, Equation (\ref{glucose}) yields $i\frac{\scalZ}{2}\phi = {\Omega^\cZ\cdot\psi}_{|_M} = (\nu\lrcorner\Omega^\cZ)\bullet\phi.$ But, 
\begin{eqnarray}
 i(\nu\lrcorner\Omega^\cZ)\bullet\phi &=&  i(\nu\cdot(\nu\lrcorner\Omega^\cZ)\cdot\psi)_{|_M}= (\nu\cdot\RicZ (\nu)\cdot\psi)_{|_M}\nonumber\\
&=& - \ricZ(\nu, \nu) \phi + \sum_{j=1}^n \ricZ(\nu, e_j) \ e_j\bullet\phi.
\label{mmm}
\end{eqnarray}
Taking the real part of the scalar product of Equation (\ref{mmm}) with $\phi$ yields $\frac{\scalZ}{2} = \ricZ(\nu, \nu)$, hence Equation (\ref{maga}) is satisfied.\\ \\
Now, let $M$ be a $\Spinc$ Riemannian manifold having a generalized Killing spinor field $\phi$ with a symmetric endomorphism $F$ on the tangent bundle $TM$. As mentioned in the introduction, it is straightforward to see that $2 T^\phi(X,Y) = -\< F(X), Y\>.$ We will study these generalized Killing spinors when the tensor $F$ is a Codazzi-Mainardi tensor, i.e., $F$ satisfies 
\begin{eqnarray}
 (\nabla^M_X F)(Y) = (\nabla^M_Y F)(X) \ \ \  \text{for}\ \ X, Y\in \Gamma(TM).
\label{aga1}
\end{eqnarray}
For this, we give the following lemma whose proof will be omitted since it is similar to Lemma $7.3$ in \cite{3}.
\begin{lemma}$\cite{3}$ Let $g_t$ be a smooth 1-parameter family of semi-Riemannian metrics on a $\Spinc$ manifold $(M^n, g = g_0)$ and let $F$ be a field of symmetric endomorphisms of $TM$. We consider the metric $g_\cZ =\<., .\> = dt^2 + g_t$ on $\cZ$ such that $g_t(X, Y) = g((\id -tF)^2 X, Y)$ for all vector fields $X, Y$ on $M$. We have $\< \RZ (U,\nu)\nu, V\>=0$ for all vector fields $U, V$ tangent to $M$ and if $F$ satisfies the Codazzi-Mainardi equation then $\< \RZ (U, V)W, \nu\> =0$ for all $U, V$ and $W$ on $\cZ$.
\label{bala}
\end{lemma}
{\bf Proof of Theorem \ref{fini}}: We define 
$\psi_{(0,x)}:=\phi_x$ via the identification 
$\Sigma_xM\cong\Sigma_{(0,x)}\cZ$ (resp. $\Sigma_{(0,x)}^+\cZ$ for $n$ odd) 
and $\psi_{(t,x)}=\tau_0^t\psi_{(0,x)}$. By Equation (\ref{zylweingarten}), the endomorphism $F$ is the Weingarten tensor of the immersion of ${\{0\}\times M}$ in $\cZ$ and hence by construction we have for all $X \in \Gamma(TM)$
\begin{equation}
\label{ano}
 \n^\SZ_X\psi|_{\{0\}\times M}= 0\, \ \ \ \ \hbox{and}\ \ \ \ \ \n^\SZ_\nu\psi\equiv 0.
\end{equation}
Since the tensor $F$ satisfies the Codazzi-Mainardi equation, Lemma \ref{bala} yields $g_\cZ (\RZ (U, V)W, \nu) =0$ for all $U, V$ and $W \in\Gamma(\cZ)$ and $g_\cZ (\RZ (X, \nu)\nu, Y) =0$ for all $X$ and $Y$ tangent to $M$. Hence  $R^\cZ(\nu,X)=0$ for all $X\in \Gamma(TM)$.
Let $X$ be a fixed 
arbitrary tangent vector field on $M$. Using (\ref{alphaa}) and (\ref{ano}) we get 
$$\n^\SZ_\nu\n^\SZ_X\psi = \mathcal{R}^{\Sigma \cZ} (\nu,X)\psi = \frac12 R^\cZ (X, \nu)\cdot\psi +\frac i2\Omega^\cZ (X,\nu)\psi =0.$$

Thus showing that the spinor field $\n^\SZ_X\psi$ is parallel along the geodesics $\RR\times \{x\}$. 
Now (\ref{ano}) shows that this spinor vanishes for $t=0$, hence it is 
zero everywhere on $\cZ$. Since $X$ is arbitrary, this shows that $\psi$ is parallel on $\cZ$.
\begin{cor}
 Let $(M^3, g)$ be a compact, oriented Riemannian manifold and $\phi$ an eigenspinor associated with the first eigenvalue $\lambda_1$ of the Dirac operator such that the Energy-Momentum tensor associated with $\phi$ is a Codazzi tensor. $M$ is a limiting manifold for (\ref{rodgy}) if and only if the generalized cylinder $\cZ^4$, equipped with the $\Spinc$ structure arising from the given one on $M$, is K\"{a}hler of positive scalar curvature and the immersion of $M$ in $\cZ$ has constant mean curvature $H$.
\end{cor}
{\bf Proof:} First, we should point out that every $3$-dimensional compact, oriented, smooth manifold has a $\Spinc$ structure. Now, if $M^3$ is a limiting manifold for (\ref{rodgy}), by Theorem \ref{fini}, the generalized cylinder has a parallel spinor whose restriction to $M$ is $\phi$. Since $\cZ$ is a $4$-dimensional $\Spinc$ manifold having parallel spinor, $\cZ$ is K\"{a}hler \cite{atiyah}. Moreover, using (\ref{glucose}), we have $$\Omega^M\bullet\phi = i\frac{\scalZ}{2} \phi  = i \frac{c_{n}}{2} |\Omega^M| \phi,$$ so $\scalZ \geqslant 0$. Finaly $H =  \frac 1n \tr(W) = \frac 1n \tr({-2 T^\phi}) = -\frac{2}{n}\lambda_1$, which is a constant. Now if the generalized cylinder is K\"{a}hler and $M$ is a compact hypersurface of constant mean curvature $H$, thus $M$ is compact hypersurface immersed in a $\Spinc$ manifold having parallel spinor with constant mean curvature. Since $\scalZ \geqslant 0$ and $\nu\lrcorner \Omega^\cZ = \RicZ (\nu) =0$, Corollary \ref{cococo} gives the result.

\end{document}